\documentclass[10pt,a4paper]{amsart}

\usepackage{amssymb}
\usepackage{a4wide}
\usepackage[matrix,arrow,curve]{xy}

\theoremstyle{plain}
\newtheorem{theorem}{Theorem}[section]
\newtheorem{lemma}[theorem]{Lemma}
\newtheorem{proposition}[theorem]{Proposition}
\newtheorem{corollary}[theorem]{Corollary}

\theoremstyle{definition}
\newtheorem{definition}[theorem]{Definition}

\let\le\leqslant
\let\ge\geqslant

\def\Z{\mathbb Z}
\def\R{\mathbb R}
\def\a{\mathbf a}
\def\c{\mathbf c}
\def\d{\mathbf d}
\def\e{\mathbf e}
\def\x{\mathbf x}
\def\y{\mathbf y}
\def\D{\mathcal D}
\def\veeL{\vee_{\!L}}

\def\Frac{{\operatorname{Frac}}}
\def\maxlen(#1){|#1|_{\max}}
\def\len{\operatorname{len}}
\def\infs{\inf{\!}_ s}
\def\sups{\sup{\!}_s}

\def\t{\operatorname{\it t}}

\begin{document}

\title{Abelian subgroups of Garside groups}

\author{Eon-Kyung Lee \and Sang Jin Lee}

\address{Department of Applied Mathematics, Sejong University,
    Seoul, Korea}
\email{eonkyung@sejong.ac.kr}

\address{Department of Mathematics, Konkuk University, Seoul, Korea}
\email{sangjin@konkuk.ac.kr}

\begin{abstract}
In this paper, we show that
for every abelian subgroup $H$ of a Garside group,
some conjugate $g^{-1}Hg$ consists of
ultra summit elements
and the centralizer of $H$ is a finite index subgroup
of the normalizer of $H$.
Combining with the results on translation numbers in Garside groups,
we obtain an easy proof of the algebraic flat torus theorem for Garside groups
and solve several algorithmic problems concerning abelian subgroups
of Garside groups.

\medskip\noindent
{\bf\em Key Words: \/}
Garside group; conjugacy class; abelian subgroup;
algebraic flat torus theorem; translation number.

\noindent
{\bf\em 2000 Mathematics Subject Classification: \/}
Primary 20F36; Secondary 20F10
\end{abstract}

\maketitle

\section{Introduction}

For mapping class groups and Artin groups of finite type,
there are several results on properties of abelian subgroups.
For mapping class groups of surfaces with negative Euler characteristic,
Birman, Lubotzky and McCarthy (1983) computed
the maximal rank of an abelian subgroup.
McCarthy (1982) showed that every abelian subgroup containing
a pseudo-Anosov mapping class is generated by two elements,
a pseudo-Anosov and a periodic mapping class.
This result plays an important role in the recent work of Birman, Gebhardt
and Gonz\'alez-Meneses (2006) on the conjugacy problem
that every pseudo-Anosov braid has a (uniformly bounded) small power
whose ultra summit set consists of rigid elements.
Recently, Hamemst\"adt (2005) and
Behrstock and Minsky (2005) solved
Brock-Farb's Rank Conjecture that the maximal rank of a quasi-flat is the same as
the maximal rank of an abelian subgroup.
For Artin groups of finite type, Charney and Peifer (2003)
showed that the maximal rank of an abelian
subgroup is equal to the number of vertices in its Coxeter graph.

However, relatively few things are known for abelian subgroups of Garside groups,
a lattice theoretic generalization of braid groups and Artin groups of finite type.
The following are what we have found in literature.

\begin{itemize}
\item
Every abelian subgroup of a Garside group is torsion-free and
finitely generated.\\
Every Garside group is torsion-free by Dehornoy (1998).
Charney, Meier and Whittlesey (2004) showed that
every Garside group has finite virtual cohomological dimension,
hence every abelian subgroup of a Garside group is finitely generated.

\item The algebraic flat torus theorem holds for Garside groups.\\
In 1995, Alonso and Bridson proved the algebraic flat torus theorem
for semihyperbolic groups:
\emph{if\/ $G$ is a semihyperbolic group and $A$ is a finitely-generated abelian group,
then every monomorphism $\phi:A\to G$ is a quasi-isometric embedding.}
Semihyperbolic groups are groups that admit a quasi-geodesic bicombing.
It is known that Garside groups are biautomatic
(Dehornoy and Paris, 1999; Dehornoy, 2002) and
biautomatic groups are semihyperbolic (Alonso and Bridson, 1995;
Bridson and Haefliger, 1999).
\end{itemize}

In this paper, we are interested in abelian subgroups of Garside groups.
For Garside groups, there are well-established theories
on solving the conjugacy problem, which involve computing
super summit sets or ultra summit sets.
Intuitively, the super summit set of an element is the set
of all conjugates that have the shortest normal form
in the conjugacy class, and the ultra summit set is a subset
of the super summit set whose elements are contained
in closed orbits under cycling.

We first show the following.

\medskip
\noindent\textbf{Theorem A} (Theorem \ref{thm:SSS} (i))\ \ \em
Let $G$ be a Garside group and $H$ an abelian subgroup of\/ $G$.
Then there exists an element $g\in G$
such that $g^{-1}Hg$ consists of ultra summit elements.
\normalfont\medskip

In particular, the group $g^{-1}Hg$ consists of super summit elements,
hence each element has the shortest normal form in the conjugacy class.
This result yields Proposition~\ref{thm:normalizer}
that for each abelian subgroup $H$ of $G$,
the centralizer $Z_G(H)$ is a finite index subgroup of
the normalizer $N_G(H)$.

Using Theorem A together with the results on
translation numbers in Garside groups,
we get an easy proof of the algebraic flat torus theorem
for Garside groups without using semihyperbolicity.

\medskip
\noindent\textbf{Theorem B} (Theorem \ref{thm:quasi-flat})\ \ \em
If\/ $G$ is a Garside group and $A$ is a finitely-generated abelian group,
then every monomorphism $\phi:A\to G$ is a quasi-isometric embedding.
\normalfont\medskip

\noindent
Furthermore, if $H$ is cyclic,
then $g^{-1}Hg$ is $(1,2)$-quasi-isometric to the real line
for some element $g\in G$
(see Proposition~\ref{prop:linear}).

\medskip
Lastly, we show that the following algorithmic problems
for abelian subgroups are solvable for Garside groups.
In the statement, a subset $\{g_1,\ldots,g_n\}$
of an abelian group is said to be linearly independent
if $g_1^{c_1}\cdots g_n^{c_n}=1$ implies $c_1=\cdots=c_n=0$,
and it is called a basis
if it forms a linearly independent set of generators.

\medskip
\noindent\textbf{Theorem C}
(Lemma~\ref{lemma:basis}--\ref{lemma:conj-prob-abel})\ \ \em
Let $G$ be a Garside group.
\begin{enumerate}
\item[(i)]
(Basis problem for abelian subgroups)\ \
There is a finite-time algorithm that,
given a collection $h_1,\ldots,h_n$
of mutually commuting elements of\/ $G$,
finds a basis of the subgroup generated by $h_1,\ldots,h_n$.

\item[(ii)]
(Membership problem for abelian subgroups)\ \
There is a finite-time algorithm that,
given an element $g$ of $G$ and a finite collection
$h_1,\ldots,h_n$ of mutually commuting elements of\/ $G$,
decides whether $g$ belongs to the subgroup generated by
$h_1,\ldots,h_n$.

\item[(iii)]
(Conjugacy membership problem for abelian subgroups)\ \
There is a finite-time algorithm that,
given an element $g$ of $G$ and a finite collection
$h_1,\ldots,h_n$ of mutually commuting elements of\/ $G$,
decides whether $g$ is conjugate to an element of the subgroup
generated by $h_1,\ldots,h_n$.

\item[(iv)]
(Equality problem for abelian subgroups)\ \
There is a finite-time algorithm that,
given two finite collections
$\{h_1,\ldots,h_n\}$ and $\{h_1',\ldots,h_m'\}$
of mutually commuting elements of\/ $G$,
decides whether they generate the same subgroup.

\item[(v)]
(Conjugacy problem for abelian subgroups)\ \
There is a finite-time algorithm that,
given two finite collections
$\{h_1,\ldots,h_n\}$ and $\{h_1',\ldots,h_m'\}$
of mutually commuting elements of\/ $G$,
decides whether they generate conjugate subgroups.
\end{enumerate}
\normalfont\medskip

We note that a simpler version of
Theorem A restricted to cyclic subgroups
was presented in our earlier preprint
titled `Stable super summit sets in Garside groups'.
Shortly after our posting, Birman informed that
she together with Gebhardt and Gonz\'alez-Meneses
independently had obtained the same result
(precisely, stable ultra summit sets in Garside groups are nonempty),
and soon posted the preprint (Birman, Gebhardt and Gonzalez-Meneses, 2006).
Our proof was more involved than theirs
because we proved the non-emptiness without using the convexity theorem
of Franco and Gonz\'alez-Meneses (2003) and Gebhardt (2005).
In this paper, we use the convexity theorem in the proof of Theorem A.

\subsection*{Acknowledgement}
This is a revised version of the paper published in
Communications in Algebra \textbf{36} (2008) 1121--1139.
The differences between the published and this versions
are in Theorem~\ref{thm:ConvexityOfSnU} and Theorem~\ref{thm:SSS}(ii).
In the published paper, we claimed that $C_U(g)$ and $C_U(H)$ are closed under
both $\wedge_L$ and $\vee_L$,
but Juan Gonz\'alez-Meneses and Volker Gebhardt
found independently that they are not closed under $\vee_L$.
The authors are very grateful to them.
We remark that this change is irrelevant to the main results---Theorems~A, B and C---of the paper.
The second author was supported by Konkuk University in 2006.

\section{Garside groups}

We start with a brief review of Garside groups.
See (Garside, 1969; Epstein et al., 1992; Birman, Ko and Lee, 1998;
Dehornoy and Paris, 1999; Dehornoy, 2002; Franco and Gonz\'alez-Meneses, 2003;
Gebhardt 2005) for details.

\subsection{Garside monoids and groups}
Let $M$ be a monoid.
Let \emph{atoms}  be the elements $a\in M\setminus \{1\}$
such that $a=bc$ implies either $b=1$ or $c=1$.
For $a\in M$, let $\maxlen(a)$ be the supremum of the lengths of all expressions of
$a$ in terms of atoms. The monoid $M$ is said to be \emph{atomic}
if it is generated by its atoms and $\maxlen(a)<\infty$ for every $a\in M$.
In an atomic monoid $M$, there are partial orders $\le_L$ and $\le_R$:
$a\le_L b$ if $ac=b$ for some $c\in M$;
$a\le_R b$ if $ca=b$ for some $c\in M$.

\begin{definition}
A finitely generated monoid $M$ is called a \emph{Garside monoid} if
\begin{enumerate}
\item[(i)] $M$ is atomic;
\item[(ii)] $M$ is left and right cancellative;
\item[(iii)] the posets $(M,\le_L)$ and $(M,\le_R)$ are lattices;
\item[(iv)] there exists an element $\Delta$, called a
\emph{Garside element}, satisfying the following:\\
(a) for each $a\in M$, $a\le_L\Delta$ if and only if $a\le_R\Delta$;\\
(b) the set $\{a\in M: a \le_L\Delta\}$ generates $M$.
\end{enumerate}
\end{definition}

An element $a\in M$ is called a \emph{simple element} if $a\le_L\Delta$.
Let $\D$ denote the set of simple elements.
Let $\wedge_L$ and $\veeL$ denote the gcd and lcm with respect to $\le_L$.

Garside monoids satisfy Ore's conditions,
and thus embed in their groups of fractions.
A \emph{Garside group} is defined as the group of fractions
of a Garside monoid.
When $M$ is a Garside monoid and $G$ the group of fractions of $M$,
we identify the elements of $M$ and their images in $G$
and call them \emph{positive elements} of $G$.
$M$ is called the \emph{positive monoid} of $G$,
often denoted $G^+$.
The partial orders $\le_L$ and $\le_R$, and thus the lattice structures
in the positive monoid $G^+$ can be extended
to the Garside group $G$ as follows:
$g\le_L h$ (resp.{} $g\le_R h$) for $g,h\in G$
if $gc=h$ (resp.{} $cg=h$) for some $c\in G^+$.

Let $\tau\colon G\to G$ be the inner automorphism of $G$ defined by
$\tau(g)=\Delta^{-1} g\Delta$.
It is known that  $\tau(G^+)=G^+$, that is, the positive monoid is
invariant under conjugation by $\Delta$.

For $g\in G$, there are integers $r\le s$ such that
$\Delta^r\le_L g\le_L\Delta^s$.
Hence, the invariants
$$
\inf(g)=\max\{r\in\Z:\Delta^r\le_L g\}
\quad\mbox{and}\quad
\sup(g)=\min\{s\in\Z:g\le_L \Delta^s\}
$$
are well-defined.
The canonical length is defined by $\len(g)=\sup(g)-\inf(g)$.
For $g\in G$, there is a unique expression
$$
g=\Delta^r s_1\cdots s_k,
$$
called the \emph{normal form} of $g$,
where $s_1,\ldots,s_k\in \D\setminus\{1,\Delta\}$ and
$(s_is_{i+1}\cdots s_k)\wedge_L \Delta=s_i$ for $i=1,\ldots,k$.
In this case, $\inf(g)=r$ and\/ $\sup(g)=r+k$.

\subsection{Conjugacy problem in Garside groups}
Let $g$ be  an element of a Garside group $G$ with
normal form $\Delta^r s_1\cdots s_k$.
The \emph{cycling} $\c(g)$ and the \emph{decycling} $\d(g)$ of $g$
are defined by
$$
\c(g)=\Delta^r s_2\cdots s_k \tau^{-r}(s_1)
\quad\mbox{and}\quad
\d(g)=\Delta^r\tau^r(s_k)s_1\cdots s_{k-1}.
$$
Let $[g]$ denote the conjugacy class of $g$ in $G$.
We define
$$
\infs(g)=\max\{\inf(h):h\in [g]\}
\quad\mbox{and}\quad
\sups(g)=\min\{\sup(h):h\in [g]\}.
$$
The \emph{super summit set} $[g]^S$ and the
\emph{ultra summit set} $[g]^U$ are defined as follows:
\begin{eqnarray*}
[g]^S&=&\{h\in [g]:\inf(h)=\infs(g),\ \sup(h)=\sups(g)\};\\{}
[g]^U&=&\{h\in [g]^S:\mbox{$\c^k(h)=h$ for some $k\ge 1$}\}.
\end{eqnarray*}
Elements of super summit sets and ultra summit sets are called {\em super summit elements} and
{\em ultra summit elements}, respectively.

\begin{lemma}\label{thm:cycl}
Let $g$ be an element of a Garside group.
\begin{enumerate}
\item[(i)] $\tau(\c(g))=\c(\tau(g))$ and $\tau(\d(g))=\d(\tau(g))$.
\item[(ii)] $\inf(g)\le\inf(\c(g))\le\sup(\c(g))\le\sup(g)$.
\item[(iii)] $\inf(g)\le\inf(\d(g))\le\sup(\d(g))\le\sup(g)$.
\item[(iv)] If\/ $\inf(g)<\infs(g)$,
  then $\inf(\c^l(g))>\inf(g)$ for some $l\ge 1$.
\item[(v)] If\/ $\sup(g)>\sups(g)$,
  then $\sup(\d^l(g))<\sup(g)$ for some $l\ge 1$.
\item[(vi)] $\c^k(\d^l(g))\in[g]^U$ for some $k,l\ge 0$.
\end{enumerate}
\end{lemma}

\begin{theorem}\label{thm:SSS-USS}
Let $g$ be an element of a Garside group.
\begin{enumerate}
\item[(i)] Both $[g]^S$ and $[g]^U$ are finite and non-empty.
\item[(ii)] Both $[g]^S$ and $[g]^U$ are closed under
$\c$, $\d$ and $\tau$.
\item[(iii)] For $h,h'\in[g]^S$,
there exists a finite sequence of elements in $[g]^S$
$$h=h_0\to h_1\to\cdots\to h_m=h'$$
such that for $i\ge1$, $h_{i}=s_i^{-1}h_{i-1}s_i$
for some $s_i\in\D$.
The same is true for $[g]^U$.
\end{enumerate}
\end{theorem}

The above theorem solves the conjugacy problem in Garside groups.
Two elements are conjugate if and only if
their super summit sets are the same
because super summit sets are non-empty by (i).
We can obtain at least one element in the super summit set by
Lemma~\ref{thm:cycl}~(vi),
and we can compute the whole super summit set from a single element by (iii).

For an element $g$ of a Garside group $G$, we define
$$
C_S(g)=\{h\in G: h^{-1}gh\in[g]^S\}\quad\mbox{and}\quad
C_U(g)=\{h\in G: h^{-1}gh\in[g]^U\}.
$$
Theorem~\ref{thm:SSS-USS}~(ii) and (iii)
come from the following theorem, known as the Convexity Theorem,
due to Franco and Gonz\'alez-Meneses (2003) for $C_S(g)$ and
Gebhardt (2005) for $C_U(g)$.
The statement is a little more general than those
in their original papers, hence we include a sketchy proof.

\begin{theorem}\label{thm:ConvexityOfSnU}
For each element $g$ of a Garside group $G$,
both $C_S(g)$ and $C_U(g)$ are closed under $\wedge_L$
and multiplication by $\Delta^{\pm 1}$ on the right.
Moreover, $C_S(g)$ is closed under $\veeL$.
\end{theorem}

\begin{proof}
It is obvious that both $C_S(g)$ and $C_U(g)$ are closed under
multiplication by $\Delta^{\pm 1}$ on the right.
Franco and Gonz\'alez-Meneses (2003) showed that
the set $C_S(g)\cap G^+$ is closed under $\wedge_L$,
and Gebhardt (2005) showed that
the set $C_U(g)\cap G^+$ is closed under $\wedge_L$.
(The original statements require that $g$ is
a super summit element and a ultra summit element, respectively,
but these conditions can be easily dropped.)

Suppose that $h_1,h_2\in C_S(g)$.
Choose an integer $u$ such that $h_1\Delta^u$ and $h_2\Delta^u$
are positive elements.
Because $C_S(g)$ is closed under multiplication by $\Delta^{\pm 1}$
on the right, both $h_1\Delta^u$ and $h_2\Delta^u$ belong to
$C_S(g)\cap G^+$.
Since
\begin{eqnarray*}
(h_1\Delta^u)\wedge_L(h_2\Delta^u)
&=&(\Delta^u\tau^u(h_1))\wedge_L(\Delta^u \tau^u(h_2))
=\Delta^u\, (\tau^u(h_1)\wedge_L\tau^u(h_2))\\
&=&\Delta^u\, (\tau^u(h_1\wedge_L h_2))
=(h_1\wedge_L h_2)\,\Delta^u
\end{eqnarray*}
belongs to $C_S(g)$, we have $h_1\wedge_L h_2\in C_S(g)$.
Therefore $C_S(g)$ is closed under $\wedge_L$.
The same arguments yield that $C_U(g)$ is closed under $\wedge_L$.

For the closedness under $\veeL$, we temporarily define
$C_S'(g)=\{h\in G: hgh^{-1}\in[g]^S\}$.
Observe that $h^{-1}gh$ in the definition of $C_S(g)$
is replaced by $hgh^{-1}$.
Using the same argument on $C_S(g)$, we can see that
$C_S'(g)$ is closed under $\wedge_R$.
Note that $h\in C_S(g)$ if and only if $h^{-1}\in C_S'(g)$.

Suppose that $h_1,h_2\in C_S(g)$.
Note that $h_1^{-1}\wedge_R h_2^{-1}=(h_1\veeL h_2)^{-1}$.
Since $h_1^{-1},h_2^{-1}\in C_S'(g)$,
$h_1^{-1}\wedge_R h_2^{-1}=(h_1\veeL h_2)^{-1}\in C_S'(g)$,
hence $h_1\veeL h_2\in C_S(g)$.
Therefore $C_S(g)$ is closed under $\veeL$.
\end{proof}

Throughout this paper, if not specified,
$G$ is always assumed to be a Garside group
with a Garside element $\Delta$, the set $\D$ of simple elements.
In addition, $L_\Delta$ denotes the maximal word length $\maxlen(\Delta)$
of the Garside element $\Delta$.

\subsection{Translation numbers}

For a finitely generated group $G$ and
a finite set $X$ of generators for $G$,
the \emph{translation number} with respect to $X$ of
an element $g\in G$ is defined by
$$\t_{G,X}(g)=\lim_{n\to \infty}\frac{|g^n|{}_X}n,$$
where $|\cdot|_X$ denotes the shortest word length
in the alphabet $X\cup X^{-1}$.
If there is no confusion about the group $G$,
we simply write $\t_X(g)$ instead of $\t_{G,X}(g)$.
The following lemma describes basic properties
of translation numbers (Gersten and Short, 1991).

\begin{lemma}\label{lemma:trans-property}
Let $G$ be a group and $X$ a finite set of generators for $G$.
\begin{itemize}
\item[(i)] For all $g\in G$, $\t_X(g)$ is well-defined.
\item[(ii)] For all $g,h\in G$, $\t_X(h^{-1}gh)=\t_X(g)$,
that is $\t_X(\cdot)$ is a conjugacy invariant.
\item[(iii)] For all $g\in G$ and $n\in\Z$,
$\t_X(g^n)=|n|\cdot \t_X(g)$.
\item[(iv)] If $g,h\in G$ commute with each other,
then $\t_X(gh)\le \t_X(g)+\t_X(h)$.
\end{itemize}
\end{lemma}

For the translation numbers in Garside groups,
the following are known (Lee, 2007; Lee and Lee, 2006a and 2006b).

\begin{theorem}\label{thm:trans-rational}
Let $G$ be a Garside group.
\begin{enumerate}
\item[(i)]
If $g$ is a super summit element of\/ $G$, then $|g|_\D -2\le \t_\D(g)\le |g|_\D$.

\item[(ii)]
The translation numbers in $G$ are
rational of the form $p/q$
for some integers $p, \ q$ with $1\le q\le L_\Delta^2$.

\item[(iii)]
If\/ $g$ is a non-identity element of\/ $G$, then $\t_\D(g)\ge 1/L_\Delta$.

\item[(iv)]
There is a finite-time algorithm that, given an element of\/ $G$,
computes its translation number.
\end{enumerate}
\end{theorem}

\section{Ultra summit property of abelian subgroups}

\begin{definition}
Two $n$-tuples $(h_1,\ldots,h_n)$ and $(h_1',\ldots,h_n')$ of elements
in a group $G$ are said to be \emph{simultaneously conjugate} if there exists
an element $g$ of $G$ such that $h_i'=g^{-1}h_ig$ for all $i=1,\ldots,n$.
Such an element $g$ is called a \emph{simultaneous conjugator} from
$(h_1,\ldots,h_n)$ to $(h_1',\ldots,h_n')$.
\end{definition}

\begin{lemma}\label{lemma:SimConj}
Let\/ $h_1,\ldots,h_n, g$ be elements of\/ $G$
such that $gh_i=h_ig$ for all $i=1,\ldots, n$.
Let $\Delta^u s_1\cdots s_k$ be the normal form of $g$.
Then the following hold.
\begin{itemize}
\item[(i)]
$(h_1,\ldots,h_n,g)$ is simultaneously conjugate to
$(h_1',\ldots,h_n',\c(g))$ by $\tau^{-u}(s_1)$
such that for each $i=1,\ldots,n$,
if\/ $h_i$ is a super/ultra summit element, then so is $h_i'$.

\item[(ii)]
$(h_1,\ldots,h_n,g)$ is simultaneously conjugate to
$(h_1',\ldots,h_n',\d(g))$ by $s_k^{-1}$
such that for each $i=1,\ldots,n$,
if\/ $h_i$ is a super/ultra summit element, then so is $h_i'$.

\item[(iii)]
$(h_1,\ldots,h_n,g)$ is simultaneously conjugate to
$(h_1',\ldots,h_n',g')$ such that $g'$ is a ultra summit element
and for each $i=1,\ldots,n$,
if\/ $h_i$ is a super/ultra summit element, then so is $h_i'$.
Furthermore, we can find a simultaneous conjugator from
$(h_1,\ldots,h_n,g)$ to $(h_1',\ldots,h_n',g')$ in finite time.
\end{itemize}
\end{lemma}

\begin{proof}
(i)\ \
Let $s=\tau^{-u}(s_1)$.
Since $s^{-1}gs=\c(g)$, $(h_1,\ldots,h_n)$ is simultaneously conjugate to
$(s^{-1}h_1s,\ldots,s^{-1}h_ns, \c(g))$.
Suppose $h_i$ is a super summit element for some $i=1,\ldots,n$.
Since $g^{-1}h_ig=h_i$, we have $g\in C_S(h_i)$.
Since $\Delta^{u+1}\in C_S(h_i)$ and $C_S(h_i)$ is closed under $\wedge_L$,
the element
$$
\Delta^{u+1}\wedge_L g=\Delta^us_1=\tau^{-u}(s_1)\Delta^u=s\Delta^u
$$
belongs to $C_S(h_i)$.
Therefore $s\in C_S(h_i)$, hence $s^{-1}h_is$ is a super summit element.
Similarly, we can show that if $h_i$ is a ultra summit element, then
so is $s^{-1}h_is$.

\medskip
(ii)\ \
It can be proved similarly to (i).

\medskip
(iii)\ \
$\c^k\d^l(g)$ belongs to the ultra summit set for some integers $k,l\ge 0$.
By (i) and (ii), $(h_1,\ldots,h_n,g)$ is simultaneously conjugate to
$(h_1',\ldots,h_n',\c^k\d^l(g))$
such that for each $i=1,\ldots,n$,
if\/ $h_i$ is a super/ultra summit element, then so is $h_i'$.
Note that the simultaneous conjugator is a product of $k+l$ elements
in $\D\cup\D^{-1}$ that can be computed from $g$.
\end{proof}

\begin{corollary}\label{cor:simultaneous}
Let $g_1,\ldots,g_n$ be mutually commuting elements of\/ $G$.
Then $(g_1,\ldots,g_n)$ is simultaneously conjugate to
$(g_1',\ldots,g_n')$ such that
each $g_i'$ is a ultra summit element.
Furthermore, we can find a simultaneous conjugator from
$(g_1,\ldots,g_n)$ to $(g_1',\ldots,g_n')$ in finite time.
\end{corollary}

\begin{proof}
Assume that $g_i$ is a ultra summit element
for all $i=1,\ldots,k$ for some $0\le k<n$.
Applying Lemma~\ref{lemma:SimConj},
we can conclude that
$(g_1,\ldots,g_k,g_{k+1},\ldots,g_n)$ is simultaneously
conjugate to $(g_1'',\ldots,g_k'',g_{k+1}'',\ldots,g_n'')$ such that
$g_i''$ is a ultra summit element for all $i=1,\ldots, k+1$.
Using induction on $k$, the desired result is obtained.
\end{proof}

Now, we generalize the notions of $C_S(\cdot )$ and $C_U(\cdot )$.
For a subset $T$ of a Garside group $G$, define
\begin{eqnarray*}
C_S(T) & = & \{x\in G: x^{-1}gx\in[g]^S \mbox{ for all } g\in T\}; \\
C_U(T) & = & \{x\in G: x^{-1}gx\in[g]^U \mbox{ for all } g\in T\}.
\end{eqnarray*}

\begin{theorem}\label{thm:SSS}
Let $H$ be an abelian subgroup of a Garside group $G$.
Then the following hold.
\begin{itemize}
\item[(i)]
There exists an element $g\in G$
such that $g^{-1}Hg$ consists of ultra summit elements.
In other words, $C_U(H)$ and hence $C_S(H)$ are nonempty.

\item[(ii)]
Both $C_S(H)$ and $C_U(H)$ are closed under $\wedge_L$
and multiplication by $\Delta^{\pm 1}$ on the right.
Moreover, $C_S(H)$ is closed under $\veeL$.
\end{itemize}
\end{theorem}

\begin{proof}
(i) \
Let $h_1,\ldots,h_n$ be a finite set of generators for $H$.
Let $\{W_1,W_2,\ldots\}$ be the set of all freely reduced words on
$\{x_1^{\pm 1},\ldots,x_n^{\pm 1}\}$ such that $W_i=x_i$ for $i=1,\ldots,n$.
In particular $H=\{W_i(h_1,\ldots,h_n):i=1,2,\ldots\}$.
For each integer $m\ge n$, let $S_m$ denote the set of all $n$-tuples
$(h_1',\ldots,h_n')$ simultaneously conjugate to $(h_1,\ldots,h_n)$
such that each $W_i(h_1',\ldots,h_n')$
is a ultra summit element for $i=1,\ldots,m$.

\medskip
Firstly, we claim that $S_m$ is nonempty for all $m\ge n$.
For $i>n$, let $h_i=W_i(h_1,\ldots,h_n)$.
Since $h_1,\dots,h_m$ mutually commute,
there exists $g\in G$ such that
$g^{-1}h_ig$ is a ultra summit element for all $i=1,\ldots,m$ by
Corollary~\ref{cor:simultaneous}.
For $1\le i\le m$, let $h_i'=g^{-1}h_ig$.
Then
$$h_i'=g^{-1}h_ig=g^{-1}W_i(h_1,\ldots,h_n)g=W_i(h_1',\ldots,h_n'),
$$
hence $(h_1',\ldots,h_n')$ belongs to $S_m$.
Therefore $S_m$ is nonempty for all $m\ge n$.

\medskip
Secondly, we claim that $S_n$ is a finite set.
If $(h_1',\ldots,h_n')\in S_n$, then each $h_i'$ is
a ultra summit element for $i=1,\ldots,n$.
Therefore $S_n$ is a subset of $[h_1]^U\times\cdots\times[h_n]^U$.
Because each ultra summit set $[h_i]^U$ is a finite set,
$S_n$ is a finite set.

\medskip
It is obvious from the definition of $S_m$
that $S_n\supset S_{n+1}\supset\cdots$.
Because $S_n$ is a finite set as observed,
$\bigcap_{m\ge n} S_m=S_{m_0}$ for some $m_0\ge n$.
Because $S_{m_0}$ is nonempty, we can take
$(h_1',\ldots,h_n')$ from $S_{m_0}$.
Since $(h_1,\ldots,h_n)$ is simultaneously
conjugate to $(h_1',\ldots,h_n')$,
there exists an element $g\in G$ such that
$h_i'=g^{-1}h_ig$ for $i=1,\ldots,n$,
hence $g^{-1}Hg=\{W_i(h_1',\ldots,h_n'):i=1,2,\ldots\}$.
Since $S_{m_0}=\bigcap_{m\ge n} S_m$,
$W_i(h_1',\ldots,h_n')$ is a ultra summit element for all $i\ge 1$,
hence $g^{-1}Hg$ consists of ultra summit elements.

\medskip\noindent
(ii) \
Since $C_S(H) = \bigcap_{h\in H} C_S(h)$ and
$C_U(H) = \bigcap_{h\in H} C_U(h)$, it is obvious due to
Theorem~\ref{thm:ConvexityOfSnU}.
\end{proof}

We apply the above result to centralizers and normalizers
of abelian subgroups.
For a subgroup $H$ of a group $G$,
let $Z_G(H)$ and $N_G(H)$ denote the centralizer and normalizer of $H$ in $G$,
that is, $Z_G(H)=\{g\in G: gh=hg\ \mbox{for all $h\in H$}\}$
and $N_G(H)=\{g\in G: g^{-1}Hg=H\}$.
It is obvious that $Z_G(H)$ is a subgroup of $N_G(H)$.

\begin{proposition}\label{thm:normalizer}
If\/ $H$ is an abelian subgroup of a Garside group $G$,
$[N_G(H):Z_G(H)]<\infty$.
\end{proposition}

\begin{proof}
By Theorem~\ref{thm:SSS},
we may assume that $H$ consists of ultra summit elements.
Let $h_1,\ldots,h_n$ be a finite set of generators for $H$.
Let
$$
H_0=\{ (h_1',\ldots,h_n') \in H^n:
\mbox{$(h_1',\ldots,h_n')$ is simultaneously
conjugate to $(h_1,\ldots,h_n)$}\}.
$$
Because $H$ consists of ultra summit elements,
$H_0$ is a subset of $[h_1]^U\times\cdots\times [h_n]^U$,
hence $H_0$ is a finite set.
Let $\rho$ denote the action of $N_G(H)$ on $H_0$ by conjugation, that is,
$$
\rho(g)(h_1',\ldots,h_n')=(g^{-1}h_1'g,\ldots,g^{-1}h_n'g)
$$
for $g\in N_G(H)$ and $(h_1',\ldots,h_n')\in H_0$.
Note that
$\rho(g)$ fixes $(h_1,\ldots,h_n)$
if and only if $g\in Z_G(H)$.
That is,
$Z_G(H)$ is the stabilizer of $(h_1,\ldots,h_n)$ in $N_G(H)$.
Consequently $[N_G(H):Z_G(H)]\le |H_0|<\infty$.
\end{proof}

Bestvina (1999) showed that if $H$ is a normal abelian subgroup
of an irreducible Artin group $G$ of finite type,
then $H$ is central. (The center of an irreducible Artin group is
infinite cyclic generated by $\Delta$ or $\Delta^2$,
hence $H$ is generated by $\Delta^k$ for some $k\in\Z$ with $\Delta^k$ central.)
He said that it answers a question of Jim Carlson, which motivated
his construction of the normal form complex and
the analysis of its geometric properties.
The result is proved in two steps:
(i) for any $h\in H$,
the conjugacy class $[h]$ is a finite set;
(ii) if an element $g\in G$ is not central,
then the conjugacy class $[g]$ is an infinite set.

It looks difficult to generalize the result of Bestvina
to normal abelian subgroups of Garside groups,
but the first statement can be strengthened as the following proposition.

\begin{proposition}
If\/ $H$ is a normal abelian subgroup of a Garside group $G$, then
$H$ consists of ultra summit elements.
Further, $[h]=[h]^U$ for every $h\in H$.
\end{proposition}

\begin{proof}
By Theorem~\ref{thm:SSS}, $g^{-1}Hg$ consists of ultra summit elements
for some $g\in G$.
Because $H$ is normal, $g^{-1}Hg=H$.
For any element $h\in H$, the conjugacy class $[h]$ is a subset of $H$
because $H$ is a normal subgroup.
\end{proof}

\section{Technical lemmas}
This section provides two technical lemmas.
They are easy to prove, but we include the proof
for completeness.

Let $R$ be $\R$ or $\Z$.
Let $V$ be an $R$-module.
Recall that a function $\Vert \cdot\Vert:V\to \R$
is called a \emph{seminorm} if
\begin{itemize}
\item $\Vert \x\Vert\ge 0$ for all $\x\in V$;
\item $\Vert r\x\Vert=|r|\cdot\Vert\x\Vert$
    for all $\x\in V$ and $r\in R$;
\item $\Vert \x+\y\Vert\le
    \Vert \x\Vert+\Vert \y\Vert$
    for all $\x,\y\in V$.
\end{itemize}
A seminorm $\Vert\cdot\Vert$ is called a \emph{norm} if $\Vert \x\Vert=0$
implies $\x=\mathbf 0$.

Let $\e_i$ denote the $i$th standard unit vector of $\R^n$.
Let $\Vert\cdot\Vert_\infty$ be the norm on $\R^n$ defined by
$$
\Vert \x\Vert_\infty=\max\{|x_1|,\ldots,|x_n|\}
\qquad \mbox{for}\ \
\x=(x_1,\ldots,x_n)\in\R^n.
$$

The following lemma is a generalization of the well-known approximation
$|x-p/q|\le 1/q^2$ of a real number $x$ by a rational number $p/q$.

\begin{lemma}\label{lemma:approx}
For any $\x\in\R^n$ and any positive integer $M$,
there exist $\a\in\Z^n$ and a positive integer $k\le M^n$ such that
$$
\Vert k\x-\a\Vert_\infty \le \frac 1M.
$$
\end{lemma}

\begin{proof}
For a real number $x$, let $\lfloor x\rfloor$ and $\Frac(x)$ denote
the integral and fractional part of $x$, that is,
$\lfloor x\rfloor$ is the largest integer less than or equal to $x$
and $\Frac(x)=x-\lfloor x\rfloor$.
For an $\R$-vector $\y=(y_1,\ldots,y_n)$,
let $\lfloor \y\rfloor$ denote the $\Z$-vector
$(\lfloor y_1\rfloor,\ldots,\lfloor y_n\rfloor)$,
and let $\Frac(\y)=\y-\lfloor\y\rfloor$.
Divide the $n$-cube $[0,1]^n$ into $M^n$ small $n$-cubes
congruent to $[0,1/M]^n$ and consider the set
$$
\{\Frac(k\x): k=0,1,\ldots,M^n\}.
$$
By the pigeonhole principle, there exist $0\le k_1<k_2\le M^n$ such that
$\Frac(k_1\x)$ and $\Frac(k_2\x)$ are contained in the same small $n$-cube,
hence
$$
\Vert\Frac(k_2\x)-\Frac(k_1\x)\Vert_\infty\le \frac1M.
$$
Let $k=k_2-k_1$ and
$\a=\lfloor k_2\x\rfloor -\lfloor k_1\x\rfloor$.
Then
$$
k\x-\a
= k_2\x- k_1\x-\lfloor k_2\x\rfloor +\lfloor k_1\x\rfloor
= \Frac(k_2\x)-\Frac(k_1\x),
$$
hence $\Vert k\x-\a\Vert_\infty \le 1/M$.
\end{proof}

It is well-known that any two norms
on a finite dimensional $\R$-vector space are equivalent,
that is, if $\Vert\cdot\Vert_1$ and $\Vert\cdot\Vert_2$ are two norms on $\R^n$,
then there exist positive constants $C_1$ and $C_2$ such that
$(1/C_1)\Vert \x\Vert_1\le \Vert \x\Vert_2\le C_2\Vert \x\Vert_1$
for all $\x\in \R^n$. Usual proofs of this inequality are not constructive,
hence they do not give the constants specifically.
Because we need to compute the constants in finite time for solving some
algorithmic problems, we show the following lemma.

\begin{lemma}\label{lemma:seminorm}
Let $\Vert\cdot\Vert:\R^n\to \R$ be a seminorm.
Suppose that there exist positive integers $K$ and $L$ such that
\begin{itemize}
\item
$\Vert\a\Vert\ge 1/L$ for all $\a\in\Z^n$ with
$1\le\Vert \a\Vert_\infty\le (2nKL)^n$;
\item
$\Vert\e_i\Vert\le K$ for all $i=1,\ldots,n$.
\end{itemize}
Let $C_1=(2L)^{n+1}(nK)^n$ and $C_2=nK$.
Then for all $\x\in\R^n$,
$$
\frac 1{C_1}\cdot\Vert\x\Vert_\infty
\le \Vert\x\Vert
\le C_2\cdot \Vert\x\Vert_\infty.
$$
In particular $\Vert\cdot\Vert$ is a norm.
\end{lemma}

\begin{proof}
Let $\x=x_1\e_1+\cdots+x_n\e_n$.
For simplicity, we may assume that $\Vert \x\Vert_\infty=1$.
Then $|x_i|\le \Vert \x\Vert_\infty=1$ for all $i$,
and the triangular inequality yields
$$
\Vert \x\Vert
=\Vert x_1\e_1+\cdots+x_n\e_n\Vert
\le |x_1|\cdot \Vert \e_1\Vert+\cdots+|x_n|\cdot \Vert \e_n\Vert
\le nK=C_2.
$$

\medskip
Applying Lemma~\ref{lemma:approx} with $M=2nKL$, we obtain
an integer $1\le k\le M^n$ and a $\Z$-vector $\a$ such that
$$
\Vert k\x-\a\Vert_\infty\le\frac1M.
$$
Since $M=2nKL\ge 2$,
\begin{eqnarray*}
\Vert \a\Vert_\infty
    &\ge& \Vert k\x\Vert_\infty - \Vert k\x-\a\Vert_\infty
    \ge k-\frac1M\ge 1-\frac1M>0;\\
\Vert\a\Vert_\infty
    &\le&\Vert k\x\Vert_\infty + \Vert k\x -\a\Vert_\infty
    \le k+\frac1M
    \le M^n+\frac1M<M^n+1.
\end{eqnarray*}
Since $\Vert\a\Vert_\infty$ is an integer,
we have $1\le\Vert\a\Vert_\infty\le M^n$,
hence by the assumption
$$
\Vert \a\Vert\ge \frac1L.
$$
Since $\Vert k\x-\a\Vert_\infty\le 1/M$,
\begin{eqnarray*}
\Vert k\x-\a\Vert
&\le&
C_2\cdot \Vert k\x-\a\Vert_\infty
\le \frac {C_2}M=\frac{nK}{2nKL}=\frac 1{2L},\\
\Vert k\x\Vert
&=& \Vert \a+(k\x-\a)\Vert
\ge \Vert \a\Vert -\Vert k\x-\a\Vert
\ge \frac 1L-\frac1{2L}=\frac1{2L},\\
\Vert\x\Vert
&=& \frac 1k\cdot\Vert k\x\Vert
\ge \frac1{M^n}\cdot \frac1{2L}
=\frac1{2LM^n}=\frac1{(2L)^{n+1}(nK)^n}=\frac1{C_1}.
\end{eqnarray*}
Now, we have proved that $1/C_1\le\Vert\x\Vert\le C_2$ for $\Vert \x\Vert_\infty=1$.
\end{proof}

\section{Quasi-flatness of abelian subgroups}

Here, we prove the algebraic flat torus theorem
for Garside groups.
Note that, in a Garside group,
translation numbers of non-identity elements
are strictly positive by Theorem~\ref{thm:trans-rational}~(iii),
hence translation numbers restricted to an abelian subgroup
give a norm by Lemma~\ref{lemma:trans-property}.

\begin{definition}
Let $G$ be a group with a finite set $S$ of generators for $G$.
For a homomorphism $\phi:\Z^n\to G$, define a seminorm
$\Vert\cdot\Vert_{\phi, S}:\Z^n\to\R$
by $\Vert\x\Vert_{\phi, S}= \t_S(\phi(\x))$ for $\x\in\Z^n$.
\end{definition}

Note that if the group $G$ is translation separable
(that is, translation numbers of non-torsion elements are
strictly positive),
then the seminorm $\Vert\cdot\Vert_{\phi, S}$ becomes a norm on $\Z^n$
if and only if $\phi$ is a monomorphism.
The following lemma implies that if a homomorphism $\phi:\Z^n\to G$ is injective
near the origin $\mathbf 0$, then $\phi$ is a monomorphism.

\begin{lemma}\label{lemma:bilip}
Let $G$ be a Garside group, and let $\phi:\Z^n\to G$ be a homomorphism.
Let $K=\max\{ \Vert \e_1\Vert_{\phi, \D},\ldots,
\Vert \e_n\Vert_{\phi, \D}\}$.
Suppose that $\phi(\a)$ is not the identity for all $\a\in\Z^n$ with
$1\le \Vert \a\Vert_\infty\le(2nKL_\Delta)^n$.
Then $\phi$ is a monomorphism and
$$
\frac1{D_1}\cdot\Vert \a\Vert_\infty
\le \Vert \a\Vert_{\phi, \D}
\le D_2\cdot \Vert \a\Vert_\infty
\qquad\mbox{for all $\a\in\Z^n$},
$$
where
$D_1=(2L_\Delta)^{n+1}(nK)^n$ and $D_2=nK$.
\end{lemma}

\begin{proof}
Recall that if $g\ne 1$, then $\t_\D(g)\ge 1/L_\Delta$
by Theorem~\ref{thm:trans-rational}.
By the assumption, if $1\le\Vert \a\Vert_\infty\le (2nKL_\Delta)^n$,
then $\Vert\a\Vert_{\phi, \D}\ge 1/L_\Delta$ since $\phi(\a)$ is not the identity.
Applying Lemma~\ref{lemma:seminorm} with $L=L_\Delta$, we have
the desired inequality.
Therefore, $\Vert\cdot\Vert_{\phi,\D}$ is a norm
and $\phi$ is a monomorphism.
\end{proof}

\begin{definition}
Let $(X,d_X)$ and $(Y,d_Y)$ be metric spaces.
A map $f:X\to Y$ is called a
\emph{$(\lambda,\epsilon)$-quasi-isometric embedding}
if there are constants $\lambda\ge 1$ and $\epsilon\ge 0$ such that
$$
\frac1\lambda\ d_X(x_1,x_2)-\epsilon
\le d_Y(f(x_1),f(x_2))
\le \lambda\ d_X(x_1,x_2)+\epsilon
\qquad
\mbox{for all $x_1,x_2\in X$}.
$$
We often suppress $(\lambda,\epsilon)$, saying just
quasi-isometric embedding.
A quasi-isometric embedding $f:X\to Y$ is called a \emph{quasi-isometry}
if there exists a constant $\delta\ge 0$ such that
each point of $Y$ is contained in the
$\delta$-neighborhood of $f(X)$.
\end{definition}

If $X$ is an $\R$- or $\Z$-module with a norm $\Vert\cdot\Vert$, then
we can define a metric on $X$ by setting
$d(x_1,x_2)=\Vert x_1-x_2\Vert$ for $x_1,x_2\in X$.
Abusing notation, $(X,\Vert\cdot\Vert)$ denotes
both the normed space and the induced metric space $(X,d)$.

\medskip

The following is the algebraic flat torus theorem
that every abelian subgroup of a Garside group
is quasi-isometric to $\Z^n$ for some $n$.

\begin{theorem}\label{thm:quasi-flat}
If\/ $G$ is a Garside group and $A$ is a finitely-generated abelian group,
then every monomorphism $\phi:A\to G$ is a quasi-isometric embedding.
\end{theorem}

\begin{proof}
Because Garside groups are torsion-free and $\phi$ is a monomorphism,
the finitely-generated abelian group $A$ is also torsion-free.
Hence we may assume $A=\Z^n$ for some integer $n$.

Let $H=\phi(A)$.
By Theorem~\ref{thm:SSS}, there exists $g\in G$
such that $g^{-1}Hg$ consists of ultra summit elements.
Let $H'=g^{-1}Hg$.
Since for every $h\in H$
$$
|h|_\D - 2|g|_\D\le |g^{-1}hg|_\D\le |h|_\D+2|g|_\D,
$$
$(H,|\cdot|_\D)$ is $(1,2|g|_\D)$-quasi-isometric to $(H',|\cdot|_\D)$.
In Theorem 7.1 of (Lee, 2007), it is shown that
if $h$ is a super summit element, then
$$
|h|_\D -2\le \t_\D(h)\le |h|_\D.
$$
Therefore, $(H',|\cdot|_\D)$ is $(1,2)$-quasi-isometric to
$(H', \t_\D(\cdot))$.
Now, we have shown that $(H,|\cdot|_\D)$ is $(1,2(|g|_\D+1))$-quasi-isometric
to $(H',\t_\D(\cdot))$.

On the other hand, $(H', \t_\D(\cdot))$ is isometric to $(\Z^n,\Vert\cdot\Vert_{\phi, \D})$
by definition of the seminorm $\Vert\cdot\Vert_{\phi, \D}$,
hence it is quasi-isometric to $(\Z^n, \Vert\cdot\Vert_\infty)$
by Lemma~\ref{lemma:bilip}.
\end{proof}

\section{Stable super summit sets}
For an element $g$ of a Garside group $G$,
we define the \emph{stable super summit set} of $g$ in $G$ as
$$
[g]^{St}=\{h\in[g]: h^n\in[g^n]^S\
\mbox{for all $n\ge 1$}\}.
$$
In other words, every power of an element in a stable super summit set
is a super summit element.
In this section, we explore elementary properties of stable super summit sets,
because stable super summit sets are useful in the study of conjugacy
classes in Garside groups
(Lee and Lee, 2006a and 2006b; Birman, Gebhardt and Gonzalez-Meneses, 2006).

For an element $g\in G$, define
$$
C_{St}(g)=\{h\in G: h^{-1}gh\in[g]^{St}\}.
$$
Applying Theorem~\ref{thm:SSS} to infinite cyclic groups,
we obtain the following theorem.

\begin{theorem}\label{thm:StSSS}
Let $g$ be an element of a Garside group $G$.
Then the following hold.
\begin{itemize}
\item[(i)]
The stable super summit set $[g]^{St}$ is nonempty.

\item[(ii)]
$C_{St}(g)$ is closed under $\wedge_L$, $\veeL$
and multiplication by $\Delta^{\pm 1}$ on the right.

\item[(iii)]
If\/ $h\in [g]^{St}$, then $\tau(h), \c(h), \d(h) \in [g]^{St}$.

\item[(iv)]
If\/ $h,h'\in[g]^{St}$, there exists a finite sequence
of elements in $[g]^{St}$
$$ h=h_0\to h_1\to\cdots\to h_m=h'$$
such that for $i=1,\ldots,m$, $h_i=s_i^{-1}h_{i-1}s_i$
for some $s_i\in\D$.
\end{itemize}
\end{theorem}

\begin{proof}
Let $H$ be the infinite cyclic group generated by $g$.
By Theorem~\ref{thm:SSS}, there exists $x\in G$ such that
$x^{-1}Hx$ consists of ultra summit elements.
Therefore $x^{-1}gx$ belongs to the stable super summit set of $g$,
hence (i) is proved.
(ii) follows from Theorem~\ref{thm:SSS}
since $C_{St}(g) = C_S(H)$.
(iii) and (iv) follow from (ii).
\end{proof}

We remark that the above theorem is not sufficient
to make a finite-time algorithm
for computing stable super summit sets,
because we need a finite-time algorithm for testing
whether an element $h\in [g]$ is contained
in the stable super summit set $[g]^{St}$:
a naive algorithm would test whether $h^n\in[g^n]^S$
for all positive integers $n$.
In (Lee and Lee, 2006b), it is shown that $h\in[g]^{St}$ if and only if
$h^n\in[g^n]^S$ for $n=1,\ldots,L_\Delta$,
in other words, the tuple $(h,h^2,h^3,\ldots,h^{L_\Delta})$ consists of super
summit elements.
Combining with Corollary~\ref{cor:simultaneous},
we obtain a finite-time algorithm for computing stable super summit sets.

It would be quite interesting to see the interplay between
the study of stable super summit sets and that of translation numbers
in Garside groups.
Non-emptiness of stable super summit sets is essential to
the study of translation numbers in (Lee and Lee, 2006a)
and of periodically geodesic powers in (Lee and Lee, 2006b),
from which a finite-time algorithm for computing stable super summit sets
comes (Lee and Lee, 2006b).

\medskip
We now estimate $\infs$ and $\sups$ of $g^{n+m}$ in terms of $\infs$ and $\sups$ of
$g^n$ and $g^m$ (in Proposition~\ref{thm:inf_estimate})
by using the following lemma.

\begin{lemma}\label{thm:ineq}
Let $h$ be an element of a Garside group $G$. For $n\ge 1$,
\begin{enumerate}
\item[(i)] $n\infs(h)\le\infs(h^n)\le n\infs(h)+(n-1)$;
\item[(ii)] $n\sups(h)-(n-1)\le\sups(h^n)\le n\sups(h)$.
\end{enumerate}
\end{lemma}

\begin{proof}
Theorem 6.1 in (Lee, 2007) states that
$$\infs(g)  \le \frac{\infs(g^n)}n < \infs(g)+1\quad\mbox{and}\quad
\sups(g)-1 <  \frac{\sups(g^n)}n \le \sups(g).
$$
Since $\infs$ and $\sups$ are integer-valued,
we get the desired inequalities.
\end{proof}

\begin{proposition}\label{thm:inf_estimate}
Let $g$ be an element of a Garside group $G$.
For $m,n\ge 1$,
\begin{enumerate}
\item[(i)] $\infs(g^m)+\infs(g^n)\le \infs(g^{m+n})\le \infs(g^m)+\infs(g^n)+1$;
\item[(ii)] $\sups(g^m)+\sups(g^n)-1\le\sups(g^{m+n})\le\sups(g^m)+\sups(g^n)$.
\end{enumerate}
\end{proposition}

\begin{proof}
We prove only (i), because (ii) can be proved similarly.

We first show that $\infs(g^{m+n})\ge \infs(g^m)+\infs(g^n)$.
Substituting $g^m$ and $g^n$ for $h$ in Lemma~\ref{thm:ineq}~(i),
$$
\infs(g^{m(m+n)})=\infs((g^m)^{m+n})\ge (m+n)\infs(g^m);
\atop
\infs(g^{n(m+n)})=\infs((g^n)^{m+n})\ge (m+n)\infs(g^n).
$$
Therefore,
\begin{equation}\label{eqn:inf-1}
\infs(g^{m(m+n)}) + \infs(g^{n(m+n)})
\ge (m+n)(\infs(g^m)+\infs(g^n)).
\end{equation}
On the other hand,
substituting $g^{m+n}$ for $h$ in Lemma~\ref{thm:ineq}~(i),
$$
\infs(g^{m(m+n)})
=\infs((g^{m+n})^m)
\le m\infs(g^{m+n})+m-1;
\atop
\infs(g^{n(m+n)})
=\infs((g^{m+n})^n)
\le n\infs(g^{m+n})+n-1.
$$
Therefore,
\begin{equation}\label{eqn:inf-2}
\infs(g^{m(m+n)})+\infs(g^{n(m+n)})
\le (m+n)(\infs(g^{m+n})+1)-2.
\end{equation}
If $\infs(g^{m+n})\le \infs(g^m)+\infs(g^n)-1$,
then (\ref{eqn:inf-2}) implies
$$
\infs(g^{m(m+n)})+\infs(g^{n(m+n)})
\le (m+n)(\infs(g^m)+\infs(g^n))-2,
$$
which contradicts (\ref{eqn:inf-1}).
Consequently, $\infs(g^{m+n})\ge \infs(g^m)+\infs(g^n)$.

\medskip

The other inequality can be proved similarly.
Substituting $g^m$ and $g^n$ for $h$ in Lemma~\ref{thm:ineq}~(i),
$$
\infs(g^{m(m+n)})
=\infs((g^m)^{m+n})
\le (m+n)\infs(g^m)+(m+n)-1;
\atop
\infs(g^{n(m+n)})
=\infs((g^n)^{m+n})
\le (m+n)\infs(g^n)+(m+n)-1.
$$
Therefore,
\begin{equation}\label{eqn:inf-3}
\infs(g^{m(m+n)}) + \infs(g^{n(m+n)})
\le (m+n)(\infs(g^m)+\infs(g^n)+2)-2.
\end{equation}
Substituting $g^{m+n}$ for $h$ in Lemma~\ref{thm:ineq}~(i),
$$
\infs(g^{m(m+n)})
=\infs((g^{m+n})^m)
\ge m\infs(g^{m+n});
\atop
\infs(g^{n(m+n)})
=\infs((g^{m+n})^n)
\ge n\infs(g^{m+n}).
$$
Therefore,
\begin{equation}\label{eqn:inf-4}
\infs(g^{m(m+n)}) + \infs(g^{n(m+n)})
\ge (m+n)\infs(g^{m+n}).
\end{equation}
If $\infs(g^{m+n})\ge \infs(g^m)+\infs(g^n)+2$,
then (\ref{eqn:inf-4}) implies
$$
\infs(g^{m(m+n)}) + \infs(g^{n(m+n)})
\ge (m+n)(\infs(g^m)+\infs(g^n)+2),
$$
which contradicts (\ref{eqn:inf-3}).
Consequently, $\infs(g^{m+n})\le \infs(g^m)+\infs(g^n)+1$.
\end{proof}

\medskip

We close this section with some remarks on stable super summit sets.

First, we remark that
every infinite cyclic subgroup of a Garside group is $(1,\epsilon)$-quasi-isometric
to the real line $\R$ for some $\epsilon\ge 0$.
Let $g$ be an element of a Garside group $G$ and let $H$ be the infinite cyclic group
generated by $g$.
If $g$ belongs to its stable super summit set,
then by Theorem~\ref{thm:trans-rational} (Theorem 7.1 of Lee, 2007)
$$
|n| \cdot\t_\D(g) \le | g^n |_\D \le |n| \cdot\t_\D(g) +2
\quad\mbox{for all $n\in\Z$}.
$$
Therefore, the infinite cyclic group $H$ is $(1,2)$-quasi-isometric
to the real line $\R$
endowed with the norm $\Vert x\Vert=\t_\D(g)\cdot |x|$ for $x\in \R$.
If $g$ does not belong to its stable super summit set, then
there exists an element $x\in G$ such that $x^{-1}gx\in [g]^{St}$.
Since $H$ and $x^{-1}Hx$ are $(1,2|x|_\D)$-quasi-isometric,
$H$ is $(1,2|x|_\D+2)$-quasi-isometric to $\R$.
Hence, we have the following proposition.

\begin{proposition}\label{prop:linear}
Let $G$ be a Garside group and $H$ an infinite cyclic subgroup of $G$.
Then there exists an element $g\in G$ such that $g^{-1}Hg$
is $(1,2)$-quasi isometric to the real line.
In particular, $H$ is $(1,\epsilon)$-quasi isometric to the real line
for some $\epsilon\ge 0$.
\end{proposition}

\begin{figure}
$$
\begin{array}{ccc}
\xymatrix{
& \sigma_1\sigma_2\sigma_3 \ar@/^.7em/[dr]^{\sigma_1}\\
  \sigma_1\sigma_3\sigma_2 \ar@/^.7em/[ur]^{\sigma_3} \ar@/_.7em/[dr]_{\sigma_1}
&&
  \sigma_2\sigma_1\sigma_3 \ar[ll]_{\sigma_2}\\
& \sigma_3\sigma_2\sigma_1 \ar@/_.7em/[ur]_{\sigma_3}
}
&\qquad&
\xymatrix{ &&\\
  \sigma_1\sigma_3\sigma_2 \ar@/_1em/[rr]^{\sigma_1\sigma_3}
&&
  \sigma_2\sigma_1\sigma_3 \ar@/_1em/[ll]_{\sigma_2}
}\\
\mbox{(a) $[g_1]^S=[g_1]^U$}
&&
\mbox{(b) $[g_1]^{St}$}
\end{array}
$$
\caption{Minimal conjugacy graphs of $[g_1]^S$, $[g_1]^U$ and $[g_1]^{St}$}
\label{fig:MinimalGraph}
\end{figure}
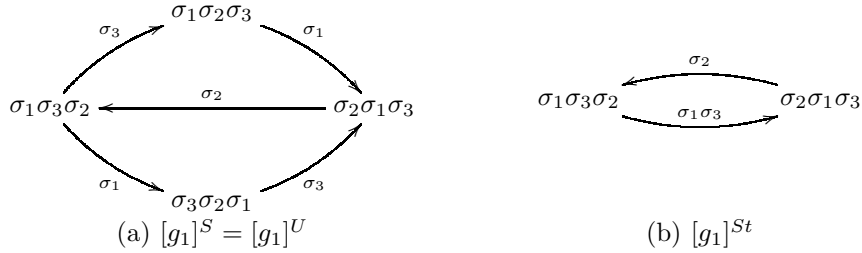

\medskip
Next, we show by an example that
(i) the stable super summit set is different from
both the super summit set and the ultra summit set;
(ii) we cannot obtain an element of the stable super summit set
by applying only cyclings and decyclings.
Consider the positive 4-braid monoid
$$ B_4^+=\langle \sigma_1,\sigma_2,\sigma_3\mid
\sigma_1\sigma_2\sigma_1=\sigma_2\sigma_1\sigma_2,\
\sigma_2\sigma_3\sigma_2=\sigma_3\sigma_2\sigma_3,\
\sigma_1\sigma_3=\sigma_3\sigma_1\rangle.
$$
This is a Garside monoid with Garside element
$\Delta=\sigma_1\sigma_2\sigma_1\sigma_3\sigma_2\sigma_1$.
Let
$g_1=\sigma_1\sigma_2\sigma_3$,
$g_2=\sigma_3\sigma_2\sigma_1$,
$g_3=\sigma_1\sigma_3\sigma_2$ and
$g_4=\sigma_2\sigma_1\sigma_3$.
Note that $g_i$'s are simple elements and conjugate to each other.
It is easy to see that
$$
[g_1]^S=[g_1]^U=\{g_1, g_2,g_3,g_4\}.
$$
Therefore, the stable super summit set of $g$ is different from
the super/ultra summit set of $g$.
The normal forms of $g_i^2$ are as follows:
$g_1^2=(\sigma_1\sigma_2\sigma_3\sigma_1\sigma_2)\sigma_3$;
$g_2^2=(\sigma_3\sigma_2\sigma_1\sigma_3\sigma_2)\sigma_1$;
$g_3^2=\Delta$; $g_4^2=\Delta$.
Therefore, $\inf(g_1^2)=\inf(g_2^2)=0$ and $\inf(g_3^2)=\inf(g_4^2)=1$.
It is easy to see that
$$
[g_1]^{St}=\{g_3,g_4\}.
$$
Note that $\c^k(g_i)=\d^k(g_i)=g_i$ for $i=1,\ldots,4$ and all $k\ge 1$.
In particular, we cannot obtain an element of the stable super summit set
by applying only cyclings and decyclings to $g_1$ or $g_2$.
Figure~\ref{fig:MinimalGraph} shows the minimal conjugacy graphs,
defined by Franco and Gonz\'alez-Meneses (2003),
of $[g_1]^S=[g_1]^U$ and $[g_1]^{St}$.

\section{Solvability of some algorithmic problems}

In this section, we prove Theorem C that several decision problems
concerning abelian subgroups are solvable for Garside groups.
Throughout this section, $G$ denotes a Garside group.

\begin{lemma}[Integer relation algorithm]
\label{lemma:integer-rel}
There is a finite-time algorithm that,
given a collection $h_1,\ldots,h_n$ of mutually commuting elements of $G$,
decides whether they are linearly independent, and if not, finds a
nontrivial word representing the identity element.
\end{lemma}

\begin{proof}
Let $\phi:\Z^n\to G$ be the homomorphism defined by $\phi(\e_i)=h_i$.
Test whether $\phi(\a)$ is the identity for $\a\in\Z^n$
with $1\le\Vert \a\Vert_\infty\le (2nKL_\Delta)^n$,
where $K=\max\{ \t_\D(h_1),\ldots, \t_\D(h_n) \}$.
If $\phi(\a)$ is not the identity for all such $\a$,
then $\phi$ is a monomorphism
by Lemma~\ref{lemma:bilip},
hence $h_1,\ldots,h_n$ are linearly independent.
If $\phi(\a)$ is the identity for some $\a$,
then $h_1,\ldots,h_n$ are linearly dependent
and $\phi(\a)$ is a nontrivial word
with word length $\le n(2nKL_\Delta)^n$
representing the identity.
\end{proof}

\begin{lemma}[Basis problem for abelian subgroups]
\label{lemma:basis}
There is a finite-time algorithm that,
given a collection $h_1,\ldots,h_n$
of mutually commuting elements of\/ $G$,
finds a basis of the subgroup generated by $h_1,\ldots,h_n$.
\end{lemma}

\begin{proof}
Let $H$ be the abelian subgroup generated by $h_1,\ldots,h_n$.
Apply Lemma~\ref{lemma:integer-rel}.
If $h_1,\ldots,h_n$ are linearly independent, we are done.
Otherwise we obtain a nontrivial word representing the identity element.
Let $\phi:\Z^n\to G$ be the homomorphism defined by $\phi(\e_i)=h_i$.
Let $\phi(\a)$ be the identity for some $\a\ne \mathbf 0$.
Because Garside groups are torsion-free, we may assume that
$\a$ is primitive, that is, the gcd of the entries of $\a$ is 1.
Then the standard algorithm using
the Hermite normal form (for example, see (Cohen, 1993))
finds $n-1$ vectors $\a_1,\ldots,\a_{n-1}$ such that
$\phi(\a_1),\ldots,\phi(\a_{n-1})$ generate $H$.
Continue the above argument to these newly obtained generators.
Induction on the number of generators completes the proof.
\end{proof}

\begin{lemma}[Membership problem for abelian subgroups]
\label{lemma:membership}
There is a finite-time algorithm that,
given an element $g$ and a finite collection
$h_1,\ldots,h_n$ of mutually commuting elements of $G$,
decides whether $g$ is contained in the subgroup generated by
$h_1,\ldots,h_n$.
\end{lemma}

\begin{proof}
Let $H$ be the subgroup generated by $h_1,\ldots,h_n$.
Because the basis problem for abelian subgroups is solvable for
Garside groups, we may assume that $\{h_1,\ldots,h_n\}$ forms a basis for $H$.
Let $\phi:\Z^n\to G$ be the homomorphism defined by $\phi(\e_i)=h_i$.
Then $\phi$ satisfies the hypothesis in Lemma~\ref{lemma:bilip}.
Compute the translation number $\t_\D(g)$.
If $g\in H$, then $\phi(\a_0)=g$ for some $\a_0\in\Z^n$.
Note that $\Vert\a_0\Vert_{\phi, \D}= \t_\D(\phi(\a_0))=\t_\D(g)$.
By Lemma~\ref{lemma:bilip},
$$
\frac1{D_2}\cdot \t_\D(g)
\le \Vert \a_0\Vert_\infty
\le D_1\cdot \t_\D(g)
$$
where $D_1$ and $D_2$ are constants
computable from $\t_\D(h_1),\ldots,\t_\D(h_n)$ and the maximal word length
$L_\Delta$ of the Garside element $\Delta$.
Therefore, to decide whether $g$ belongs to $H$, it suffices to
test whether $\phi(\a)=g$ for $\a\in\Z^n$ with
$\frac1{D_2}\cdot \t_\D(g)\le \Vert \a\Vert_\infty
\le D_1\cdot \t_\D(g)$.
Because there are only finitely many such $\a$'s and
the word problem is solvable for Garside groups, we are done.
\end{proof}

\begin{lemma}[Conjugacy membership problem for abelian subgroups]
\label{lemma:membership-conj}
There is a finite-time algorithm that,
given an element $g$ and a finite collection
$h_1,\ldots,h_n$ of mutually commuting elements of $G$,
decides whether $g$ is conjugate to an element in the subgroup generated by
$h_1,\ldots,h_n$.
\end{lemma}

\begin{proof}
The same proof as Lemma~\ref{lemma:membership} solves the conjugacy membership
problem, because the translation number is conjugacy invariant
and the conjugacy problem is solvable for Garside groups.
\end{proof}

\begin{lemma}[Equality problem for abelian subgroups]
\label{lemma:equality}
There is a finite-time algorithm that,
given two finite collections $\{h_1,\ldots,h_n\}$ and $\{h_1',\ldots,h_m'\}$
of mutually commuting elements of $G$,
decides whether they generate the same subgroup.
\end{lemma}

\begin{proof}
Let $H$ and $H'$ be the abelian subgroups generated by
$\{h_1,\ldots,h_n\}$ and $\{h_1',\ldots,h_m'\}$,
respectively.
Because the membership problem is solvable,
we can decide whether $h_i$ is contained in $H'$ for $i=1,\ldots,n$.
Therefore, we can decide whether $H$ is a subgroup of $H'$.
Similarly, we can also decide whether $H'$ is a subgroup of $H$.
\end{proof}

In the following lemma, we will use the known fact that,
the simultaneous conjugacy problem is solvable for Garside groups,
that is, there is a finite-time algorithm that,
given two $n$-tuples $(h_1,\ldots,h_n)$ and $(h_1',\ldots,h_n')$
of elements in a Garside group,
decides whether they are simultaneously conjugate,
and finds a simultaneous conjugator if so.
When each tuple consists of mutually commuting elements,
the simultaneous conjugacy problem is solvable
by an algorithm similar to that for the ordinary conjugacy problem:
we can transform (by Corollay~\ref{cor:simultaneous})
each tuple into another tuple which is simultaneously conjugate
to the original one and consists of ultra summit elements,
and then use the convexity theorem.
For general case, see (Lee and Lee, 2002; Gonz\'alezez-Meneses, 2005).

\begin{lemma}[Conjugacy problem for abelian subgroups]
\label{lemma:conj-prob-abel}
There is a finite-time algorithm that,
given two finite collections
$\{h_1,\ldots,h_n\}$ and $\{h_1',\ldots,h_m'\}$
of mutually commuting elements of $G$,
decides whether they generate conjugate subgroups.
\end{lemma}

\begin{proof}
Let $H$ and $H'$ be the abelian subgroups generated by
$\{h_1,\ldots,h_n\}$ and $\{h_1',\ldots,h_m'\}$, respectively.
Since the basis problem is solvable, we may assume that
$\{h_1,\ldots,h_n\}$ and $\{h_1',\ldots,h_m'\}$ are
bases of $H$ and $H'$.
If $n\neq m$, it is clear that $H$ and $H'$ are not conjugate.
Therefore, we may assume that $n=m$.

Let
$K=\max\{ \t_\D(h_1),\ldots, \t_\D(h_n)\}$,
$K'=\max\{ \t_\D(h_1'),\ldots, \t_\D(h_n')\}$
and $D_1 = (2L_\Delta)^{n+1}(nK)^n$.
Let
$$
T=\{ h_1^{k_1}\cdots h_n^{k_n}: |k_i|\le D_1K'
\ \mbox{for $i=1,\ldots,n$}\}.
$$
Applying Lemma~\ref{lemma:bilip} to the monomorphism
$\phi:\Z^n\to H\subset G$ defined by $\phi(\e_i)=h_i$,
we can see that if $h\in H$ and $\t_\D(h)\le K'$, then $h\in T$.
Since $\t_\D(\cdot)$ is a conjugacy invariant,
the subset $T$ contains the union $\cup_{i=1}^n ([h_i']\cap H)$ of
the sets of conjugates of $h_i'$ in $H$.

Note that the subgroups $H$ and $H'$ are conjugate if and only if
there exists an $n$-tuple $(h_i'',\ldots,h_n'')$ of elements in $H$ such that
\begin{itemize}
\item[(i)] $(h_i'',\ldots,h_n'')$ is simultaneously conjugate to $(h_1',\ldots,h_n')$ and
\item[(ii)] $\{h_i'',\ldots,h_n''\}$ forms a basis for $H$.
\end{itemize}
The $n$-tuple $(h_i'',\ldots,h_n'')\in H^n$ satisfying the property (i)
belongs to $T^n$ as observed.
Since both the simultaneous conjugacy problem and the basis problem are solvable
for Garside groups,
we can check,
for each element $(h_1'',\ldots,h_n'')$ of $T^n$,
whether the properties (i) and (ii) hold,
in a finite number of steps.
Since $T^n$ is a finite set, we are done.
\end{proof}

\end{document}